\documentclass[journal]{IEEEtran}
\usepackage{cite}

\ifCLASSINFOpdf
\else
\fi
\usepackage{array}
\usepackage{amssymb,amsmath}
\usepackage{xcolor,graphicx,float}
\usepackage{amsthm}
\usepackage{epstopdf}
\allowdisplaybreaks[4]
\usepackage{amsmath}
\usepackage{flushend}

\emergencystretch=\maxdimen
\newtheorem{theorem}{Theorem}
\newtheorem{definition}{Definition}

\newtheorem{lemma}{Lemma}
\newtheorem{assumption}{Assumption}
\newtheorem{remark}{Remark}

\newcommand{\col}{$\upshape{col}$}
\newcommand{\blk}{$\upshape{blk}$}
\newcommand{\diag}{$\upshape{diag}$}

\begin{document}

\title{Distributed strategy-updating rules for aggregative games of multi-integrator systems with coupled constraints}

\author{Xin~Cai,
        Feng~Xiao,~\IEEEmembership{Member,~IEEE,}
        and~Bo~Wei,% <-this % stops a space

\thanks{This work was supported by the National Natural Science Foundation of China (NSFC, Grant Nos. 61873074, 61903140). Corresponding author: F. Xiao.}
\thanks{X. Cai and F. Xiao are with the State Key Laboratory of Alternate Electrical Power System with Renewable Energy Sources and with the School of Control and Computer Engineering, North China Electric Power University, Beijing 102206, China (email: caixin\_xd@126.com; fengxiao@ncepu.edu.cn). X. Cai is also with the School of Electrical Engineering, Xinjiang University, Urumqi 830047, China.}

\thanks{B. Wei is with the School of Control and Computer Engineering,
North China Electric Power University,
Beijing 102206, China (email: bowei@ncepu.edu.cn).}
}

%\markboth{Transactions on Automatic Control}%
%{Cai \MakeLowercase{\textit{et al.}}: Distributed strategy-updating rules for aggregative games of multi-integrator systems with coupled constraints}
%% The only time the second header will appear is for the odd numbered pages
% *** Note that you probably will NOT want to include the author's %***
% *** name in the headers of peer review papers.                   %***

% If you want to put a publisher's ID mark on the page you can do %it like
% Remember, if you use this you must call \IEEEpubidadjcol in the %second


\IEEEtitleabstractindextext{%
\begin{abstract}
In this paper, we explore aggregative games over networks
of multi-integrator agents with coupled constraints. To reach the
general Nash equilibrium of an aggregative game, a distributed strategy-updating
rule is proposed by a combination of the coordination of
Lagrange multipliers and the estimation of the aggregator. Each player
has only access to partial-decision information and communicates with
his neighbors in a weight-balanced digraph which characterizes
players' preferences as to the values of information received from
neighbors. We first consider networks of double-integrator agents
and then focus on multi-integrator agents. The effectiveness of the
proposed strategy-updating rules is demonstrated by analyzing the convergence of corresponding dynamical systems via the Lyapunov stability
theory, singular perturbation theory and passive theory. Numerical
examples are given to illustrate our results.
\end{abstract}

}

\maketitle

\IEEEdisplaynontitleabstractindextext

% \fi
\IEEEpeerreviewmaketitle

\section{Introduction}

Distributed decision making in networked control systems modeled in the framework of game theory has attracted an increasing interest in various fields, such as economic markets \cite{Frihauf.2012}, sensor networks \cite{Stankovic.2012}, communication networks \cite{Wang.2014}, mechanical systems \cite{Deng.2019} and smart grids \cite{Gharesifard.2016}.

As a specific class of noncooperative games, aggregative games have been widely applied in engineering scenarios, due to the favorable games' property that the decision of each player is affected by some aggregation of all players' decisions. A number of distributed algorithms were proposed to seek Nash equilibrium of aggregative games. Some algorithms based on best response dynamics were applied in the demand-response scheme of smart grids \cite{Parise.2020} and the spectrum sharing in communication networks \cite{Zhou.2017}. The other algorithms based on gradient dynamics were applied to the networked Cournot competition \cite{Koshal.2016} and power allocation of small cell networks \cite{Shokri.2020}. Furthermore, there have been reports about aggregative games of high-order dynamical systems recently. Deng and Liang designed
distributed algorithms for the coordination of heterogeneous Euler-Lagrange systems \cite{Deng.2019}. Zhang \textit{et al.} studied aggregative games of nonlinear dynamic systems and devised a distributed algorithm to seek Nash equilibrium \cite{Zhang.2019}.

In various networked scenarios, physical constraints have been widely considered due to the limited capability of individuals and the capacity constraints of overall networks. However, the above mentioned literature does not consider the games with coupled constraints. In this paper, we focus on the aggregative games with coupled constraints, i.e. the feasible strategy set of each player depends on the strategies of other players. The Nash equilibrium of games in this case is referred to as generalized Nash equilibrium (GNE).

The existing distributed or decentralized GNE seeking algorithms (strategy-updating rules) for aggregative games with coupled constraints can be categorized into discrete-time and continuous-time settings depending on the properties of the systems. For discrete-time algorithms, the fixed-point iteration \cite{Grammatico.2017} is usually applied to approximate the equilibrium in a decentralized framework which requires a coordinator to collect some global information and send it to each agent. Besides, projection dynamics \cite{Belgioioso.2018} and subgradient \cite{Belgioioso.2017} are used to deal with local constraints and nonsmooth cost functions, respectively. Recently, by means of consensus estimators, distributed discrete-time algorithms based on local communications were proposed in \cite{Parise.2019,Belgioioso.2019}. Instead, the continuous-time GNE seeking algorithms, which we focus on in this paper, studied recently are mainly distributed and can be considered in control systems. Thus, the methods in control theory can be applied to the analysis and design of seeking algorithms \cite{Persis.2019a}. Deng and Nian designed a distributed algorithm for aggregative games with the coupled equality constraint, which modeled the Nash-Cournot game of generation systems \cite{Deng.2018}. Taken coupled and private constraints into consideration, an algorithm was designed by local communications with relative information for the demand response management \cite{Liang.2017}.

However, the common feature of all algorithms in aforementioned literature about games with coupled constraints is that the dynamics of players were usually assumed as first-order systems or were not considered in the games. Nonetheless,
in many engineering scenarios, the behaviors of agents are not only driven by their inherent dynamics, but also by their interests which may be in conflict with each other. For example, demand/supply response management of distributed  energy resources (DERs) was studied in \cite{Gharesifard.2016}. Given a request of demand/supply, each DER makes its decision according to its interest which depends on its decision and the price. It is supposed that DERs are price-anticipating, i.e., the pricing function depends on the average amount of energy that DERs in networks decide to consume or produce. Then, the decision process of DERs can be modeled by an aggregative game. It is known that generation systems as DERs are a class of complex dynamic systems. In this case, how to design algorithms for generation systems to maximize their interests is an interesting topic. Inspired by this context, Persis and Monshizadeh proposed an algorithm to steer a network of second-order systems to a predefined Cournot-Nash equilibrium \cite{Persis.2019c}. These observations motivate us to study the aggregative games of multi-integrator systems with coupled constraint.

This paper is to investigate aggregative games of networks of multi-integrator agents with both coupled and local constraints. In such a case, we will design a continuous-time distributed strategy-updating rule for each player and prove that the strategies of all players will reach the GNE of aggregative games with coupled equality constraints. To the best of our knowledge, this is the first paper that proposes distributed algorithms for aggregative games of networked multi-integrator agents with coupled constraints. The interactions
among agents are expressed in two terms. One is the coupled term of
estimations of Lagrange multipliers and an aggregator, and the other is the gradient term of cost functions. To deal with local constraints, different from the projected dynamics studied in \cite{Liang.2017,Deng.2019b,Persis.2019b,Lu.2019} which are not suitable for modeling agents with more complex dynamics, projected output feedbacks are utilized to overcome the technical difficulty. Then, the estimation of the aggregator and agents' dynamics are designed in two time-scales to track the aggregator quickly for stable regulation of strategies, which is different from gain regulation in estimators \cite{Liang.2017}. Moreover, based on primal-dual theory and partially coupled constraint information, the Lagrange multiplier of each agent is regulated by communicating with neighbors to make its strategy
satisfy the coupled constraints. The main contributions of this paper are summarized as follows.

1) We formulate aggregative games for networks of multi-integrator dynamics with coupled equality constraints.

2) To handle both local and coupled constraints, we propose
a strategy-updating rule including the projected output feedback,
the coordination of Lagrange multiplier based on primal-dual
theory, and the estimation of the aggregator. The design of two time-scales ensures that the strategies of all players converge exponentially to the GNE of aggregative games. Different from the previous studies based on undirected graphs,  weight-balanced digraphs are used to characterize the belief of players in the information received from their neighbors.

3) Unlike the existing literature that adopts Lyapunov stability theory to analyze the convergence of generalized Nash equilibrium seeking algorithms, this paper synthesizes singular perturbation theory, passive theory and Lyapunov stability theory to analyze the stability of the closed-loop system.

%3) All players only need to share the estimations of aggregator and Lagrange multiplier with neighbors, but not to share their strategies and costs. A key feature is that the private information of all players are protected by the proposed rules.

This paper is organized as follows. In Section \uppercase\expandafter{\romannumeral2}, the considered problem is formulated. In Section \uppercase\expandafter{\romannumeral3}, a strategy-updating rule for double-integrator agents is designed and analyzed. The rule is extended for multi-integrator agents in Section \uppercase\expandafter{\romannumeral4}. Section \uppercase\expandafter{\romannumeral5} provides simulation examples. Finally, the conclusions and future topics are stated in Section \uppercase\expandafter{\romannumeral6}.

Notations: $\mathbb{R}$ denotes the real numbers set. $\mathbb{R}^n$ is the $n$-dimensional Euclidean space. Given vector $x \in\mathbb{R}^n$, $\|x\|$ is the Euclidean norm. $\otimes$ denotes the Kronecker product. $A^T$ and $\|A\|$ are the transpose and the spectral norm of matrix $A$, respectively.  $\lambda_i(A)$ is the $i$th eigenvalue of matrix $A$. $\lambda_{2}(\cdot)$ and $\lambda_{min}(\cdot)$ are the second smallest and minimal eigenvalues, respectively.  $\col(x_1,\ldots,x_n)=[x_1^T,\ldots,x_n^T]^T$. Given matrices $A_1,\ldots,A_n$, $\blk\{A_1,\ldots,A_n\}$ denotes the block diagonal matrix with $A_i$ on the diagonal. $I_n$ is the $n\times n$ identity matrix. $\boldsymbol{1}_n$ and $\boldsymbol{0}_n$ are the $n$-dimensional column vectors with entries being ones and zeros, respectively. $\mathbf 0$ denotes a matrix consisting of all zeros with an appropriate dimension.

\section{Problem formulation}

\subsection{Aggregative games}

Consider an aggregative game $G=(\mathcal{I},\Omega,J)$ with $N$ players indexed by the set $\mathcal{I}=\{1,\ldots,N\}$. $\Omega=\Omega_1\times\cdots\times\Omega_N\subset\mathbb{R}^{Nn}$ is the strategy space of the game, where $\Omega_i\subset \mathbb{R}^{n}$ is the strategy set of player $i\in \mathcal{I}$. $J$ $=$ $(J_1,\ldots,J_N)$, where $J_i(y_i,\sigma(y))$ $:$ $\Omega_i\times \mathbb{R}^{m}\rightarrow \mathbb{R}$ is the cost function of player $i$ depending on his strategy $y_i$ $\in$ $\Omega_i$ and the aggregation of all players' strategy $\sigma(y)$. Here, $y=(y_i,y_{-i})=(y_1,\ldots,y_{i-1},y_i,y_{i+1},\ldots,y_N)\in \Omega$ denotes the strategy profile of all players, and $y_{-i}=\col(y_1,\ldots,y_{i-1},y_{i+1},\ldots,y_N)$ is a vector including strategies of all players except player $i$. $\sigma(\cdot): \Omega \rightarrow \mathbb{R}^{m}$ denotes the aggregator defined by $\sigma(y)=\sum_{i=1}^N \varphi_i(y_i)$, where $\varphi_i:\mathbb{R}^{n}\rightarrow \mathbb{R}^{m}$ is a (nonlinear) map for the local contribution of player $i$ to the aggregator. All players communicate with each other in a directed graph (digraph) $\mathcal{G}$. For a strongly connected and weight-balanced digraph $\mathcal{G}$, the Laplacian matrix $L$ has the following properties \cite[Theorem 1.37]{Bullo.2009}:  (i) $\boldsymbol{1}_N^T L=\boldsymbol{0}_N^T$; (ii) $L+L^T$ is positive semidefinite. For the detailed concepts related to graphs, please refer to \cite{Godsil.2001}.

We assume that each agent (player) in the network can be modeled by the following multi-integrator system with order $r>1$,
\begin{equation} \label{ad}
\left\{\begin{array}{l}
x_i^{(r)}=u_i,\\
y_i=P_{\Omega_i}(x_i),
\end{array}\right.
\end{equation}
where $x_i^{(r)}=\frac{d^{r}x_i}{dt^r}$. $x_i\in\mathbb{R}^{n}$, $u_i\in\mathbb{R}^{n}$, and $y_i\in\Omega_i$ are the state, control input, and output of agent $i$, respectively.

Each agent in the game regards the output $y_i$ as its strategy which belongs to a compact strategy set. Then, the projection operation $P_{\Omega_i}(\cdot): \mathbb{R}^n\rightarrow \Omega_i$ is utilized to ensure that the output is always contained in the strategy set. Moreover, the strategies of all players are coupled by the following set of linear constraints:
\begin{equation} \label{lc}
C=\{y\in \mathbb{R}^{Nn}\mid \sum_{i=1}^N A_iy_i=\sum_{i=1}^N d_i\}=\{y\in \mathbb{R}^{Nn}\mid Ay=d\},
\end{equation}
where $A_i\in\mathbb{R}^{l\times n}$, $d_i \in\mathbb{R}^{l}$, $A=[A_1,\ldots,A_N]$, $d=\sum_{i=1}^N d_i$. In Cournot competition, aggregate function $\sigma(y)$ describes the pricing function depending on decisions of all players. Thus, the net cost function is denoted by $J_i(y_i,\sigma(y))$. The balance between supply and demand can be characterized by the coupled constraint \eqref{lc}, where $d_i$ represents the demand satisfied by agent $i$. In summary, the aim of player $i$ is to choose his strategy minimizing his cost function $J_i$ and satisfying the linear coupled constraints \eqref{lc} simultaneously. The optimization problem faced by player $i$ can be described as follows.
\begin{equation} \label{op}
\begin{array}{l}
\min_{y_i\in \Omega_i} J_i(y_i,\sigma(y))\\
s.t. \ Ay=d.
\end{array}
\end{equation}

The objective in this paper is to design a strategy-updating rule for each player with dynamics \eqref{ad} to achieve the GNE of the aggregative game $G=(\mathcal{I},\Omega,J)$. To implement the distributed setting, each player updates his strategy with the local information from communicating with neighbors on a digraph $\mathcal{G}$.

The definition of GNE is given as follows and it is a natural extension of the Nash equilibrium.

\begin{definition}
\cite[Definition 3.7]{Basar.1999}\cite{Francisco.2007} For an aggregative game $G=(\mathcal{I},\Omega,J)$, a strategy profile $y^*=(y_i^*,y_{-i}^*)$ is called the generalized Nash equilibrium of the game if
\begin{equation*}
J_i(y_i^*,\sigma(y_i^*,y_{-i}^*))\leq J_i(y_i,\sigma(y_i,y_{-i}^*))
\end{equation*}
holds for any $y_i\in \Omega_i$ with $(y_i,y_{-i}^*)\in C$ and all $i\in \mathcal{I}$.
\end{definition}

At the GNE point, no player would like to change his strategy unilaterally for less cost. Note that only a finite number of players are considered in this paper. In the case of $N\rightarrow\infty$, the effect of change in $y_i$ on $\sigma(y)$ can be negligible, the equilibrium of the game is called Wardrop equilibrium \cite{Paccagnan.2019}, which is beyond the purpose of this paper.

Denote $J_i(y)=J_i(y_i,\sigma(y)), \forall i\in \mathcal{I}$, for convenience. The gradient of cost function $J_i(y)$ with respect to $y_i$ is defined by $\nabla_{y_i}J_i(y)$. Denote the pseudo-gradient mapping $F(y)=\col(\nabla_{y_1}J_1(y),\ldots,\nabla_{y_N}J_N(y))$. Some basic assumptions on the cost function $J_i$ and the pseudo-gradient $F(y)$, which are widely used in \cite{Deng.2019,Zhang.2019,Liang.2017,Persis.2019b,Lu.2019}, are respectively given as follows.

\begin{assumption} \label{as1}
For all $i\in \mathcal{I}$, $\Omega_i \subseteq\mathbb{R}^n$ is a nonempty, closed and convex set. The cost function $J_i(y_i,\sigma(y))$ is continuously differentiable and convex in $y_i$ for every fixed $y_{-i}$. And the feasible strategy set $\Omega\cap C$ is nonempty.
\end{assumption}

\begin{assumption} \label{as2}
The pseudo-gradient $F(y):\Omega\rightarrow \mathbb{R}^{Nn}$ is strongly monotone with constant $w>0$ and Lipschitz continuous with Lipschitz constant $\theta>0$.
\end{assumption}

\begin{assumption} \label{as3}
The digraph $\mathcal{G}$ is strongly connected and weight balanced.
\end{assumption}

\begin{remark}
Here, we consider the weight-balanced digraph instead of the undirected graph. The weights on the communication topology characterize how players evaluate the information received from their neighbors \cite{Nedic.2009}, according to how much they trust their neighbors. The belief of each player in their neighbors may be different from each other according to the personal preference. Compared with undirected graphs in which the common weight reflects identical belief of players in their neighbors, weight-balanced digraphs can express more characteristics about players.
\end{remark}

\subsection{GNE to Variational Inequalities}

The GNE problem of games can be reformulated as a variational inequality problem \cite{Francisco.2007}. Given a closed and convex set $\Omega$ and the mapping $F(y):\Omega\rightarrow \mathbb{R}^{Nn}$, the variational inequality, denoted by $\text{VI}(\Omega,F)$, is to find a vector $y^*\in\Omega$ such that
\begin{equation*}
(y'-y^*)^TF(y^*)\geq 0, \forall y'\in\Omega,
\end{equation*}
and the solutions to $\text{VI}(\Omega,F)$ are GNE of games, which also called variational equilibria \cite[Theorem 2.1]{Francisco.2007}. The set of solutions of VI is denoted by $\text{SOL}(\Omega,F)$.
The solution $y^*$ of $\text{VI}(\Omega,F)$ can be reformulated based on the fixed point theorem as follows:
\begin{equation}
y^* \in \text{SOL}(\Omega,F)\Leftrightarrow y^*=P_\Omega(y^*-F(y^*)).
\end{equation}

\begin{lemma}\label{lemma1}
\cite[Corollary 2.2.5 and Theorem 2.3.3]{Francisco.2003} Given the $\text{VI}(\Omega,F)$, $\Omega\subset \mathbb{R}^{Nn}$ is a convex set and the mapping $F:\Omega\rightarrow \mathbb{R}^{Nn}$ is continuous. The following statements related to solutions hold:

(1) if $\Omega$ is closed, then $\text{SOL}(\Omega,F)$ is nonempty and compact;

(2) if $\Omega$ is compact and $F(x)$ is strongly monotone, then $\text{VI}(\Omega,F)$ has a unique solution.
\end{lemma}

Here we give some results about the GNE of aggregative games, which can be derived from the results given in \cite{Lu.2019} without considering the local inequality constraints.

\begin{lemma} \label{lemma2}
Under Assumptions \ref{as1} and \ref{as2}, for the problem \eqref{op}, there exist a unique GNE $y^*=(y_i^*,y_{-i}^*)$ and a common Lagrange multiplier $\mu^*\in \mathbb{R}^l$ such that
\begin{subequations}
\begin{align}
y^*&=P_\Omega(y^*-F(y^*)-A^T\mu^*), \label{kkt1}\\
Ay^*&=d \label{kkt2}.
\end{align}
\end{subequations}
\end{lemma}

\textbf{Proof}:
For the player $i$, given strategies of the opponents, the problem \eqref{op} is an equality constrained minimization problem. Suppose $y^*$ be a GNE of game $G=(\mathcal{I},\Omega,J)$. If a suitable constraint qualification holds, there is a Lagrange multiplier $\mu_i^*\in \mathbb{R}^l$ such that the following KKT conditions are satisfied \cite[Section 5.5]{Boyd.2004}.
\begin{equation} \label{cop}
\begin{split}
\nabla_{y_i}J_i(y_i^*,y_{-i}^*)+A_i^T\mu_i^*&=0,\\
Ay^*&=d.
\end{split}
\end{equation}

For the game, due to the continuity of $F(y)$ and closed and convex set $C$, the problem of GNE can be equivalent to the VI problem. Suppose $\bar{y}$ is a solution of $\text{VI}(\Omega\cap C,F)$. If a suitable constraint qualification holds, there exists a Lagrange multiplier $\mu^* \in \mathbb{R}^l$ such that
\begin{equation} \label{vikkt}
\begin{split}
F(\bar{y})+A^T\mu^*&=0,\\
A\bar{y}&=d.
\end{split}
\end{equation}

From \eqref{cop} and \eqref{vikkt}, it follows that $\bar{y}=y^*$ if and only if $\mu_1^*=\cdots=\mu_N^*=\mu^*$ \cite[Theorem 3.1]{Francisco.2007}. Recall that the definition of $\Omega$ and $C$, under Assumption \ref{as2}, the solution of $\text{VI}(\Omega\cap C,F)$ is unique by Lemma \ref{lemma1}.  Thus, \eqref{kkt2} is satisfied. For the convex strategy set of each player, the solution of $\text{VI}(\Omega\cap C,F)$ can be characterized by the Fixed-point Theorem, which yields \eqref{kkt1}.
$\hfill\blacksquare$

\begin{remark}
The variational equilibrium is a refinement of Nash equilibria of the game\cite{Kulkarni.2012}. Although, there exist Nash equilibria for $\mu_i^*\neq \mu_j^*$, $\forall i,j\in \mathcal{I}$. The variational equilibrium is ``more socially stable'' than other equilibria of the game \cite{Francisco.2007}. The GNE in the following sections refers to the variational equilibrium.
\end{remark}

\section{GNE seeking for double-integrator agents}
In this section, we consider the agents with double-integrator dynamics, i.e., $r=2$ in \eqref{ad},which can model mobile robots, UVAs, Euler-Lagrange systems and so on.

The designed strategy-updating rule mainly has three parts including the strategy update \eqref{sup1}, the coordination of Lagrange multiplier \eqref{lme} and the aggregator estimation \eqref{dac}. The former two parts are slow systems and the later one is a fast system, and they imply that the designed rule is executed in two time-scales.  Let $\mu_i\in\mathbb{R}^l$ be the Lagrange multiplier, let $\eta_i\in\mathbb{R}^m$ be agent $i$'s estimation of the aggregator $\sigma(y)$ of agent $i$, let $z_i$ and $w_i$ be auxiliary variables, and let $k_i$ and $\alpha$ be positive constants to be designed. The strategy-updating rule for player $i$, $i\in \mathcal{I}$, is designed as follows:
\begin{equation} \label{sup1}
\begin{aligned}
\dot{x}_i&=v_i,\\
\dot{v}_i&=-k_iv_i-x_i+y_i-\nabla_{y_i}J_i(y_i,\eta_i)-A_i^T\mu_i,\\
y_i&=P_{\Omega_i}(x_i),
\end{aligned}
\end{equation}
where the cost function $J_i(y_i,\eta_i)$ indicates the outcome of the game between player $i$ and the estimation of the aggregator, $x_i(0)\in \Omega_i$, and $v_i(0)\in \mathbb{R}^n$. The coordination of the Lagrange multiplier $\mu_i$ associated with the coupled constraints is given by
\begin{equation} \label{lme}
\begin{aligned}
\dot{\mu}_i&=-\alpha\sum_{j=1}^Na_{ij}(\mu_i-\mu_j)-z_i+A_iy_i-d_i,\\
\dot{z}_i&=\alpha\sum_{j=1}^Na_{ij}(\mu_i-\mu_j),
\end{aligned}
\end{equation}
where $\mu_i$ and $z_i$ start from $\mu_i(0)=\boldsymbol{0}_l$ and $z_i(0)\in \mathbb{R}^l$ with $\sum_{i=1}^Nz_i(0)=\boldsymbol{0}_l$, respectively.
Based on the dynamic average consensus protocol in \cite{Freeman.2006,Ye.2017}, the estimation of aggregator $\sigma$ is expressed as
\begin{equation} \label{dac}
\begin{split}
\varepsilon\dot{\eta}_i&=-\eta_i-\sum_{j=1}^Na_{ij}(\eta_i-\eta_j)-\sum_{j=1}^Na_{ij}(w_i-w_j)+N\varphi_i(y_i),\\
\varepsilon\dot{w}_i&=\sum_{j=1}^Na_{ij}(\eta_i-\eta_j),
\end{split}
\end{equation}
where $\varepsilon$ is a small positive constant. The limitation of the protocol \eqref{dac} is that the tracked signal $\varphi_i$ is required to change slowly or to be constant. In order to estimate the aggregator effectively, $\varepsilon$ is needed to make the protocol be a fast system, so that  $\varphi_i$ changes relatively slowly. Note that, the estimations of Lagrange multipliers and the aggregator can be realized by the embedded technology. It is feasible to design $\alpha$ and $\varepsilon$ for the strategies of all agents to reach the GNE of game $G=(\mathcal{I},\Omega,J)$.

Let $x$ $=$ $\col(x_1,\ldots,x_N)$, $v$ $=$ $\col(v_1,\ldots,v_N)$, $y$ $=$ $\col(y_1,\ldots,y_N)$, $\eta=$ $\col(\eta_1,\ldots,\eta_N)$, $\boldsymbol{\mu}=$ $\col(\mu_1,\ldots,\mu_N)$, $\nabla J(y,\eta)=$ $ \col(\nabla_{y_1}J_1(y_1,\eta_1),\ldots,\nabla_{y_N}J_N(y_N,\eta_N))$, $z$ $=$ $\col(z_1,\ldots,z_N)$, $w$ $=$ $\col(w_1,\ldots,w_N)$, $\varphi(y)$ $=$ $\col(\varphi_1(y_1),\ldots,\varphi_N(y_N))$, $k=\diag\{k_1,\ldots,k_N\}$, $\mathbf{A}$ $=$ $\blk\{A_1,\ldots,A_N\}$, and $D$ $=$ $\col(d_1,\ldots,d_N)$. Then, \eqref{sup1}, \eqref{lme}, and \eqref{dac} can be written as
 \begin{equation} \label{cf1}
 \begin{aligned}
 \dot{x}&=v,\\
 \dot{v}&=-(k\otimes I_n)v-x+y-\nabla J(y,\eta)-\mathbf{A}^T\boldsymbol{\mu},\\
 y&=P_\Omega(x),\\
 \dot{\boldsymbol{\mu}}&=-\alpha (L\otimes I_l)\boldsymbol{\mu}-z+\mathbf{A}y-D,\\
 \dot{z}&=\alpha (L\otimes I_l)\boldsymbol{\mu},\\
 \varepsilon \dot{\eta}&=-\eta-(L\otimes I_m)\eta-(L\otimes I_m)w+N\varphi(y),\\
  \varepsilon \dot{w}&=(L\otimes I_m)\eta.
 \end{aligned}
 \end{equation}

Next, by Lemma \ref{lemma2}, the relationship between the equilibrium point of system \eqref{cf1} and the GNE of aggregative game $G=(\mathcal{I},\Omega,J)$ is analyzed. We have the following result.

\begin{lemma} \label{lemma3}
Under Assumptions \ref{as1}-\ref{as3}, $y^*$ is a GNE of aggregative game $G=(\mathcal{I},\Omega,J)$ if and only if there exist $x^*\in\mathbb{R}^{Nn}$, $\mu^*\in\mathbb{R}^l$, $z^*\in\mathbb{R}^{Nl}$, $\eta^*\in\mathbb{R}^{Nm}$, and $w^*\in\mathbb{R}^{Nm}$ such that $(x^*,\boldsymbol{0}_{Nn},\boldsymbol{1}_N\otimes\mu^*,z^*,\eta^*,w^*)$ is an equilibrium of system \eqref{cf1}.
\end{lemma}

\textbf{Proof}:
\emph{Sufficiency}: Considering the equilibrium point of system \eqref{cf1}, we have
\begin{subequations}
\begin{align}
\boldsymbol{0}_{Nn}&=v^*, \label{eq11}\\
\boldsymbol{0}_{Nn}&=-(k\otimes I_n)v^*-x^*+y^*-\nabla J(y^*,\eta^*)-\mathbf{A}^T\boldsymbol{\mu}^*,\label{eq12}\\
y^*&=P_\Omega(x^*),\label{eq13}\\
\boldsymbol{0}_{Nl}&=-\alpha (L\otimes I_l)\boldsymbol{\mu}^*-z^*+\mathbf{A}y^*-D,\label{eq14}\\
\boldsymbol{0}_{Nl}&=\alpha (L\otimes I_l)\boldsymbol{\mu}^*,\label{eq15}\\
\boldsymbol{0}_{Nm}&=-\eta^*-(L\otimes I_m)\eta^*-(L\otimes I_m)w^*+N\varphi(y^*),\label{eq16}\\
\boldsymbol{0}_{Nm}&=(L\otimes I_m)\eta^*.\label{eq17}
\end{align}
\end{subequations}

From the bottom up, under Assumption \ref{as3}, \eqref{eq17} indicates $\eta_1^*=\cdots=\eta_N^*$. Left-multiply both sides of \eqref{eq16} by $(\boldsymbol{1}_N^T\otimes I_m)$, it derives that $\eta_1^*=\cdots=\eta_N^*=\sum_{i=1}^N \varphi_i(y_i^*)$. Thus, we have $\nabla J(y^*,\eta^*)=\nabla J(y^*,\sigma(y^*))$, which implies $\nabla J(y^*,\eta^*)=F(y^*)$. Then, \eqref{eq15} indicates
$\mu_1^*=\cdots=\mu_N^*=\mu^*$. According to $\dot{z}=\alpha (L\otimes I_l)\boldsymbol{\mu}$ and $\boldsymbol{1}_N^TL=0$ for a weight-balanced digraph, we have $\sum_{i=1}^N\dot{z}_i=\boldsymbol{0}_l$, which implies $\sum_{i=1}^Nz_i(t)=\sum_{i=1}^Nz_i(0)=\boldsymbol{0}_l$, $\forall t\geq0$. Left-multiply both sides of \eqref{eq14} by $(\boldsymbol{1}_N^T\otimes I_l)$, it derives that $\sum_{i=1}^N A_iy_i^*=\sum_{i=1}^N d_i$, which satisfies \eqref{kkt2}. In addition, it follows from \eqref{eq11}, \eqref{eq12}, and \eqref{eq13} that $y^*=P_\Omega(y^*-F(y^*)-A^T\mu^*)$, which satisfies \eqref{kkt1}. Therefore, by Lemma \ref{lemma2}, it indicates that $y^*$ is the GNE of game $G=(\mathcal{I},\Omega,J)$.

\emph{Necessary}: Suppose that $y^*$ is a GNE of aggregative game $G=(\mathcal{I},\Omega,J)$. According to condition \eqref{kkt1} in Lemma \ref{lemma2}, we have $F(y^*)=\nabla J(y^*,\eta^*)$, which implies that $\eta_1^*=\cdots=\eta_N^*=\sigma(y^*)=\sum_{i=1}^N \varphi_i(y_i^*)$ satisfying \eqref{eq16}, i.e., the estimation of aggregator is identical to the true one. And it also implies that $\mu_1^*=\cdots=\mu_N^*=\mu^*$, which satisfies \eqref{eq15}. Thus, there exist some $w^*$ and $z^*$ that satisfy \eqref{eq16} and \eqref{eq14}, respectively. Meanwhile, $v^*=\boldsymbol{0}_{Nn}$ indicates that \eqref{eq12} is satisfied. So, $(x^*,\boldsymbol{0}_{Nn},\boldsymbol{1}_N\otimes\mu^*,z^*,\eta^*,w^*)$ is the equilibrium point of system \eqref{cf1}.  $\hfill\blacksquare$

Lemma \ref{lemma3} reveals that the GNE of game $G=(\mathcal{I},\Omega,J)$ can be obtained by the strategy-updating rule \eqref{sup1} if the states of system \eqref{cf1} can converge to the equilibrium point. Accordingly, we analyze the convergence of system \eqref{cf1} to illustrate the effectiveness of strategy-updating rule \eqref{sup1}. Let $\underline{k}=\min\{k_1,\ldots,k_N\}$ and $\overline{k}=\max\{k_1,\ldots,k_N\}$. We have the following conclusion.

\begin{theorem} \label{theorem1}
Suppose that Assumptions 1-3 hold and the parameters $\alpha$, $\underline{k}$ and $\overline{k}$ satisfy the following conditions:
\begin{equation*}
\begin{aligned}
\overline{k}&<3\underline{k},\ \
\underline{k}>\frac{a_1+1}{\sqrt{6a_1/5}},\\
\|\mathbf{A}\|^2&<\underline{k}(2\omega-\theta^2)-2\overline{k},\ \
\alpha>\frac{\underline{k}\|\mathbf{A}\|^2+2}{\lambda_2}.
\end{aligned}
\end{equation*}
All the agents with dynamics (1) follow the strategy-updating rule (8)-(10). Then, there exists a positive constant $\varepsilon^*$ such that for each $\varepsilon \in (0,\varepsilon^*)$, $\eta(t)$ exponentially converges to $\sum_{i=1}^N \varphi_i(y_i)\boldsymbol{1}_N$, all players' strategies $y$ globally exponentially converge to the GNE $y^*$, and the Lagrange multipliers converge to the common value $\mu^*$.
\end{theorem}

\textbf{Proof}:
First, the equilibrium point of system \eqref{cf1} is transferred to the origin. Denote $\tilde{x}=x-x^*$, $\tilde{v}=v-v^*$, $\tilde{\boldsymbol{\mu}}=\boldsymbol{\mu}-\boldsymbol{\mu}^*$, $\tilde{z}=z-z^*$, $\tilde{y}=y-y^*$. Then, system \eqref{cf1} can be rewritten as
\begin{equation}\label{cfo1}
\begin{aligned}
\dot{\tilde{x}}&=\tilde{v},\\
\dot{\tilde{v}}&=-(k\otimes I_n)\tilde{v}-\tilde{x}+\tilde{y}-(\nabla J(y,\eta)-\nabla J(y^*,\eta^*))-\mathbf{A}^T\tilde{\boldsymbol{\mu}},\\
\dot{\tilde{\boldsymbol{\mu}}}&=-\alpha (L\otimes I_l)\tilde{\boldsymbol{\mu}}-\tilde{z}+\mathbf{A}\tilde{y},\\
\dot{\tilde{z}}&=\alpha (L\otimes I_l)\tilde{\boldsymbol{\mu}},\\
\varepsilon\dot{\eta}&=-\eta-(L\otimes I_m)\eta-(L\otimes I_m)w+N\varphi(y),\\
\varepsilon\dot{w}&=(L\otimes I_m)\eta.
\end{aligned}
\end{equation}

Based on the analysis method of singular perturbations, let $\varepsilon=0$ freeze $\eta$ and $\varphi(y)$, respectively, in the reduced and boundary-layer systems, which are analyzed as follows.

1) Quasi-steady state analysis: Define $\bar{\eta}_i$ and $\bar{w}_i$ ($i\in\mathcal{I}$) as quasi-steady states. $\eta$ and $w$ are frozen at the quasi-steady states with $\eta_i=\bar{\eta}_i=\sum_{i=1}^N \varphi_i(y_i)$ and $w_i=\bar{w}_i$ by $\varepsilon=0$. In this case, the reduced system is
\begin{equation} \label{rs1}
\begin{aligned}
\dot{\tilde{x}}&=\tilde{v},\\
\dot{\tilde{v}}&=-(k\otimes I_n)\tilde{v}-\tilde{x}+\tilde{y}-(F(y)-F(y^*))-\mathbf{A}^T\tilde{\boldsymbol{\mu}},\\
\dot{\tilde{\boldsymbol{\mu}}}&=-\alpha (L\otimes I_l)\tilde{\boldsymbol{\mu}}-\tilde{z}+\mathbf{A}\tilde{y},\\
\dot{\tilde{z}}&=\alpha (L\otimes I_l)\tilde{\boldsymbol{\mu}}.\\
\end{aligned}
\end{equation}

Consider a Lyapunov function as follows:
\begin{equation*}
V=\frac{1}{2}(a_1\|\tilde{v}\|^2+\|(k\otimes I_n)\tilde{x}+\tilde{v}\|^2+\|\tilde{\boldsymbol{\mu}}\|^2+\|\tilde{\boldsymbol{\mu}}+\tilde{z}\|^2),
\end{equation*}
where $a_1>0$.

The derivative of $V$ along the trajectories of the reduced system \eqref{rs1} is
\begin{align*}
\dot{V}&=-a_1\tilde{v}^T(k\otimes I_n)\tilde{v}-a_1\tilde{v}^T\tilde{x}+a_1\tilde{v}^T\tilde{y}\\
&\ \ \ -a_1\tilde{v}^T(F(y)-F(y^*))-a_1\tilde{v}^T\mathbf{A}^T\tilde{\boldsymbol{\mu}}-\tilde{x}^T(k\otimes I_n)\tilde{x}\\
&\ \ \ +\tilde{x}^T(k\otimes I_n)\tilde{y}-\tilde{x}^T(k\otimes I_n)(F(y)-F(y^*))\\
&\ \ \ -\tilde{x}^T(k\otimes I_n)\mathbf{A}^T\tilde{\boldsymbol{\mu}}-\tilde{v}^T\tilde{x}+\tilde{v}^T\tilde{y}-\tilde{v}^T(F(y)-F(y^*))\\
&\ \ \ -\tilde{v}^T\mathbf{A}^T\tilde{\boldsymbol{\mu}}-\frac{1}{2}(\alpha\tilde{\boldsymbol{\mu}}^T(L^T+L)\otimes I_l\tilde{\boldsymbol{\mu}}-\tilde{z}^T\tilde{\boldsymbol{\mu}}+\tilde{\boldsymbol{\mu}}^T\mathbf{A}\tilde{y})\\
&\ \ \ -\tilde{\boldsymbol{\mu}}^T\tilde{z}+\tilde{\boldsymbol{\mu}}\mathbf{A}\tilde{y}-\tilde{z}^T\tilde{z}+\tilde{z}^T\mathbf{A}\tilde{y}.
\end{align*}

Let $\tilde{\boldsymbol{\mu}}^{\parallel}=\frac{1}{N}\boldsymbol{1}_N\boldsymbol{1}^T_N\otimes I_l\tilde{\boldsymbol{\mu}}$ and
$\tilde{\boldsymbol{\mu}}^{\perp}=(I_N-\frac{1}{N}\boldsymbol{1}_N\boldsymbol{1}^T_N)\otimes I_l\tilde{\boldsymbol{\mu}}$.
Then, $\tilde{\boldsymbol{\mu}}\in \mathbb{R}^{Nl}$ can be decomposed as $\tilde{\boldsymbol{\mu}}=\tilde{\boldsymbol{\mu}}^{\parallel}+\tilde{\boldsymbol{\mu}}^{\perp}$. Thus, $\tilde{\boldsymbol{\mu}}^{\parallel}=\boldsymbol{1}_N\otimes \mu$, for some $\mu \in \mathbb{R}^l$, so that $\frac{1}{2}(L+L^T)\otimes I_l\tilde{\boldsymbol{\mu}}^{\parallel}=\boldsymbol{0}_N$, and $(\frac{1}{2}\tilde{\boldsymbol{\mu}}^{\perp})^T(L+L^T)\otimes I_l\tilde{\boldsymbol{\mu}}^{\perp}\geq \lambda_2\|\tilde{\boldsymbol{\mu}}^{\perp}\|^2$,
where $\lambda_2$ is the second smallest eigenvalue of symmetric Laplacian matrix $\frac{1}{2}(L+L^T)$.

Since $\boldsymbol{1}_N^T L=\boldsymbol{0}_N$, it follows from \eqref{rs1} that $\dot{\tilde{\boldsymbol{\mu}}}^{\parallel}=\boldsymbol{0}$. If $\tilde{\boldsymbol{\mu}}^{\parallel}(0)=\boldsymbol{0}$, we have that $\tilde{\boldsymbol{\mu}}^{\parallel}(t)=0$, and $\tilde{\boldsymbol{\mu}}=\tilde{\boldsymbol{\mu}}^{\perp}$ for $t\geq0$. Thus,
\begin{align*}
\dot{V}&\leq-\underline{k}\|\tilde{x}\|^2-\underline{k}a_1\|\tilde{v}\|^2-\alpha \lambda_2\|\tilde{\boldsymbol{\mu}}\|^2-\|\tilde{z}\|^2\\
&\ \ \ -(a_1+1)\tilde{v}^T(F(y)-F(y^*))-\underline{k}\tilde{x}^T(F(y)-F(y^*))\\
&\ \ \ -\frac{(a_1+1)}{2}\tilde{v}^T\tilde{x}+(a_1+1)\tilde{v}^T\tilde{y}+\overline{k}\tilde{x}^T\tilde{y}\\
&\ \ \ -\frac{(a_1+1)}{2}\tilde{v}^T\tilde{x}-(a_1+1)\tilde{v}^T\mathbf{A}^T\tilde{\boldsymbol{\mu}}-\underline{k}\tilde{x}^T\mathbf{A}^T\tilde{\boldsymbol{\mu}}\\
&\ \ \ -2\tilde{\boldsymbol{\mu}}^T\tilde{z}+2\tilde{\boldsymbol{\mu}}^T\mathbf{A}\tilde{y}+\tilde{z}^T\mathbf{A}\tilde{y},
\end{align*}

By Assumption \ref{as2}, $\|F(y)-F(y^*)\|\leq\theta\|y-y^*\|=\theta\|\tilde{y}\|$.
%\begin{align*}
%\dot{V}&\leq-\underline{k}\|\tilde{x}\|^2-\underline{k}a_1\|\tilde{v}\|^2-\alpha \hat{\lambda}_2\|\tilde{\boldsymbol{\lambda}}\|^2-\|\tilde{z}\|^2\\
%&\ \ \ +\theta(a_1+1)\|\tilde{v}\| \|\tilde{y}\|-\underline{k}\tilde{x}^T(F(y)-F(y^*))\\
%&\ \ \ -\frac{(a_1+1)}{2}\tilde{v}^T\tilde{x}+(a_1+1)\tilde{v}^T\tilde{y}+\overline{k}\tilde{x}^T\tilde{y}\\
%&\ \ \ -\frac{(a_1+1)}{2}\tilde{v}^T\tilde{x}-(a_1+1)\tilde{v}^T\mathbf{A}^T\tilde{\boldsymbol{\lambda}}-\underline{k}\tilde{x}^T\mathbf{A}^T\tilde{\boldsymbol{\lambda}}\\
%&\ \ \ -2\tilde{\boldsymbol{\lambda}}^T\tilde{z}+2\tilde{\boldsymbol{\lambda}}^T\mathbf{A}\tilde{y}+\tilde{z}^T\mathbf{A}\tilde{y}.
%\end{align*}
And, by the strong monotonicity of $F(y)$ stated in Assumption \ref{as2} and property of the projection operator $(P_{\Omega}(x)-P_{\Omega}(x^*))^T(x-x^*)\geq\|P_{\Omega}(x)-P_{\Omega}(x^*)\|^2$, we have that $(x-x^*)^T(F(y)-F(y^*))\geq \omega\|y-y^*\|^2$. It follows that
\begin{align*}
\dot{V}&\leq-\underline{k}\|\tilde{x}\|^2-\underline{k}a_1\|\tilde{v}\|^2-\alpha \lambda_2\|\tilde{\boldsymbol{\mu}}\|^2-\|\tilde{z}\|^2\\
&\ \ \ -\underline{k}\omega\|\tilde{y}\|^2+\theta(a_1+1)\|\tilde{v}\| \|\tilde{y}\|\\
&\ \ \ -\frac{(a_1+1)}{2}\tilde{v}^T\tilde{x}+(a_1+1)\tilde{v}^T\tilde{y}+\overline{k}\tilde{x}^T\tilde{y}\\
&\ \ \ -\frac{(a_1+1)}{2}\tilde{v}^T\tilde{x}-(a_1+1)\tilde{v}^T\mathbf{A}^T\tilde{\boldsymbol{\mu}}-\underline{k}\tilde{x}^T\mathbf{A}^T\tilde{\boldsymbol{\mu}}\\
&\ \ \ -2\tilde{\boldsymbol{\mu}}^T\tilde{z}+2\tilde{\boldsymbol{\mu}}^T\mathbf{A}\tilde{y}+\tilde{z}^T\mathbf{A}\tilde{y}\\
&=-(\frac{3\underline{k}}{4}-\frac{\overline{k}}{4})\|\tilde{x}\|^2-(\underline{k}a_1-\frac{(a_1+1)^2}{4\underline{k}}-\frac{(a_1+1)^2}{4\overline{k}})\|\tilde{v}\|^2\\
&\ \ \ -\frac{1}{2}\|\tilde{z}\|^2-(\alpha \lambda_2-\underline{k}\|\mathbf{A}\|^2-2)\|\tilde{\boldsymbol{\mu}}\|^2\\
&\ \ \ -(\underline{k}\omega-\overline{k}-\frac{\|\mathbf{A}\|^2}{2})\|\tilde{y}\|^2+\theta(a_1+1)\|\tilde{v}\|\|\tilde{y}\|\\
&\ \ \ -\underline{k}\|\frac{a_1+1}{2\underline{k}}\tilde{v}+\frac{1}{2}\tilde{x}+\mathbf{A}^T\tilde{\boldsymbol{\mu}}\|^2-\overline{k}\|\frac{a_1+1}{2\overline{k}}\tilde{v}+\frac{1}{2}\tilde{x}-\tilde{y}\|^2\\
&\ \ \ -\|\sqrt{2}\tilde{\boldsymbol{\mu}}+\frac{\sqrt{2}}{2}\tilde{z}-\frac{\sqrt{2}}{2}\mathbf{A}\tilde{y}\|^2.
\end{align*}

Using Young's Inequality, it yields that
\begin{align*}
\theta(a_1+1)\|\tilde{v}\|\|\tilde{y}\|\leq\frac{(a_1+1)^2}{2\underline{k}}\|\tilde{v}\|^2+\frac{\underline{k}\theta^2}{2}\|\tilde{y}\|^2.
\end{align*}

 Then,
\begin{align*}
\dot{V}&\leq-(\frac{3\underline{k}}{4}-\frac{\overline{k}}{4})\|\tilde{x}\|^2-(\underline{k}a_1-\frac{3(a_1+1)^2}{4\underline{k}}-\frac{(a_1+1)^2}{4\overline{k}})\|\tilde{v}\|^2\\
&\ \ \ -\frac{1}{2}\|\tilde{z}\|^2-(\alpha \lambda_2-\underline{k}\|\mathbf{A}\|^2-2)\|\tilde{\boldsymbol{\mu}}\|^2\\
&\ \ \ -(\underline{k}\omega-\overline{k}-\frac{\underline{k}\theta^2}{2}-\frac{\|\mathbf{A}\|^2}{2})-\frac{\|\mathbf{A}\|^2}{2})\|\tilde{y}\|^2.
\end{align*}
The sufficient condition for $\dot{V}<0$ is that $\overline{k}<3\underline{k}$, $\underline{k}>\frac{a_1+1}{\sqrt{6a_1/5}}$,
$\|\mathbf{A}\|^2<\underline{k}(2\omega-\theta^2)-2\overline{k}$, and
$\alpha>\frac{2\underline{k}\|\mathbf{A}\|^2+4}{2\lambda_2}$. Thus, $(\tilde{x},\tilde{v},\tilde{\boldsymbol{\mu}},\tilde{z})$ globally exponentially converges to the origin.

2) Boundary-layer analysis: The boundary-layer system of \eqref{cfo1} is described in $\tau$-time scale by $\tau=t/\varepsilon$.
\begin{equation} \label{bl}
\begin{bmatrix}
\frac{d\eta}{d\tau}\\
\frac{dw}{d\tau}
\end{bmatrix}=
\begin{bmatrix}
-I_{Nm}-L\otimes I_m & -L\otimes I_m\\
L\otimes I_m & \boldsymbol{0}
\end{bmatrix}
\begin{bmatrix}
\eta\\
w
\end{bmatrix}
+
\begin{bmatrix}
N\varphi(y)\\
\boldsymbol{0}
\end{bmatrix}.
\end{equation}
According to Theorem 5 in \cite{Freeman.2006}, if the digraph $\mathcal{G}$ is strongly connected and weighted-balanced, it is clear that $\eta$ converges exponentially to $\boldsymbol{1}_N\otimes\sum_{i=1}^N \varphi_i(y_i)$.

Therefore, by Theorem 11.4 in \cite{Khalil.2002}, we have that there exists a positive constant $\varepsilon^*$ such that for all $\varepsilon \in(0,\varepsilon^*)$, the decisions of all players following the strategy-updating rule \eqref{sup1} exponentially converge to the GNE of the aggregative game.  $\hfill\blacksquare$

\section{GNE seeking for multi-integrator agents}
Taking more complex dynamics of agents (e.g., generation systems in DERs and higher order game dynamics \cite{laraki.2013}) into account, in this section, we will extend the results in the above section to the case with multi-integrator agents ($r>2$) and apply the passive theory to analyze the proposed strategy-updating rule. For ease to expound, let $n=1$ and $x_i=[x_{i1},\ldots,x_{ir}]^T\in \mathbb{R}^r$, $\forall i\in\mathcal{I}$.

The dynamics of multi-integrator agents can be described as the following linear system
\begin{equation} \label{ad2}
\begin{aligned}
\dot{x}_i&=\bar{A}x_i+\bar{B}u_i,\\
y_i&=P_{\Omega_i}(\bar{C}x_i),
\end{aligned}
\end{equation}
where $\bar{A}=[\begin{smallmatrix}\boldsymbol{0} & I_{r-1}\\ 0 & \boldsymbol{0}\end{smallmatrix}]$, $\bar{B}=[0,\ldots,0,1]^T$, and $\bar{C}=[1,0,\ldots,0]$.

Based on the strategy-updating rule \eqref{sup1} designed in Section IV, the strategy-updating rule for multi-integrator agents is given by
\begin{equation} \label{sup2}
\begin{aligned}
\dot{x}_i&=\bar{A}x_i+\bar{B}(-K_ix_i+y_i-\nabla_{y_i} J_i(y_i,\eta_i)-A_i^T\lambda_i),\\
y_i&=P_{\Omega_i}(\bar{C}x_i),
\end{aligned}
\end{equation}
where $K_i=[1,k_{i1},\ldots,k_{i(r-1)}]$ is the state feedback matrix with $k_{ij}>1, j\in\{1,\ldots,r-1\}$  to ensure that transfer function matrix $G(s)=\bar{C}(sI-(\bar{A}-\bar{B}K_i))^{-1}\bar{B}$ is strictly positive real, i.e., poles of all elements of $G(s)$ have negative real parts.

The coordination of Lagrange multiplier $\mu_i$ and estimation of aggregator $\sigma$ refer to \eqref{lme} and \eqref{dac}, respectively. Similar to the analysis in the previous section, let $H_i=\bar{A}-\bar{B}K_i$ and $H=\blk\{H_1,\ldots,H_N\}$, \eqref{sup2}, \eqref{lme}, and \eqref{dac} can be described as
\begin{equation} \label{cf2}
\begin{aligned}
\dot{x}&=Hx+(\bar{B}\otimes I_N)(y-\nabla J(y,\eta)-\mathbf{A}^T\boldsymbol{\mu}),\\
y&=P_{\Omega}((\bar{C}\otimes I_N)x),\\
\dot{\boldsymbol{\mu}}&=-\alpha (L\otimes I_l)\boldsymbol{\mu}-z+\mathbf{A}y-D,\\
 \dot{z}&=\alpha (L\otimes I_l)\boldsymbol{\mu},\\
 \varepsilon \dot{\eta}&=-\eta-(L\otimes I_m)\eta-(L\otimes I_m)w+N\varphi(y),\\
  \varepsilon \dot{w}&=(L\otimes I_m)\eta.
\end{aligned}
\end{equation}
 Since the relationship between the equilibrium point of system \eqref{cf2} and the GNE of aggregative game $G=(\mathcal{I},\Omega,J)$ is similar to what is stated in Lemma \ref{lemma3}, we do not repeat it here. The conclusion about the convergence of states in system \eqref{cf2} is given as follows.

\begin{theorem}
Suppose that Assumptions 1-3 hold and one of the following two conditions are satisfied.

1) For $\omega\geq\frac{\|\mathbf{A}\|^2}{2}$, $
\lambda_{min}(P)>3,
\alpha>\frac{4+\|\mathbf{A}\|^2}{2\lambda_2}.
$

2) For $\omega<\frac{\|\mathbf{A}\|^2}{2}$, $
\lambda_{min}(P)>3+2\|\mathbf{A}\|^2,
\alpha>\frac{4+\|\mathbf{A}\|^2}{2\lambda_2}.$
\\
All the agents with dynamics (16) follow the strategy-updating rule (17), (9) and (10). Then, there exists a positive constant $\varepsilon^*$ such that for each $\varepsilon\in(0,\varepsilon^*)$, $\eta(t)$ exponentially converges to $\sum_{i=1}^N \varphi_i(y_i)\boldsymbol{1}_N$, all players' strategies $y$ globally exponentially converge to the GNE $y^*$, and the Lagrange multipliers converge to the common value $\mu^*$.
\end{theorem}

\textbf{Proof}:
System \eqref{cf2} is a singular perturbation system with parameter $\varepsilon$.
Let $(x^*,\boldsymbol{\mu}^*,z^*)$ be the equilibrium point of the system \eqref{cf2}, and $(\bar{\eta},\bar{w})$ be the quasi-steady state. Denote $\tilde{x}=x-x^*$, $\tilde{y}=y-y^*$, $\tilde{\boldsymbol{\mu}}=\boldsymbol{\mu}-\boldsymbol{\mu}^*$, and $\tilde{z}=z-z^*$. Similar to the proof of Theorem \ref{theorem1}, the reduced system is
%\begin{equation}
%\begin{aligned}
%\dot{\tilde{x}}&=H\tilde{x}+(\bar{B}\otimes I_N)(y-y^*\\
%&\ \ \ -(\nabla J(y,\eta)-\nabla J(y^*,\eta^*))-\mathbf{A}^T\tilde{\boldsymbol{\mu}}),\\
%\tilde{y}&=P_\Omega((\bar{C}\otimes I_N)\tilde{x}),\\
%\dot{\tilde{\boldsymbol{\mu}}}&=-\alpha (L\otimes I_l)\tilde{\boldsymbol{\mu}}-\tilde{z}+\mathbf{A}(y-y^*),\\
%\dot{\tilde{z}}&=\alpha (L\otimes I_l)\tilde{\boldsymbol{\mu}},\\
%\varepsilon\dot{\eta}&=-\eta-(L\otimes I_m)\eta-(L\otimes I_m)w+N\varphi(y),\\
%\varepsilon\dot{w}&=(L\otimes I_m)\eta.
%\end{aligned}
%\end{equation}
%
%The following analysis is similar to that shown in the proof of Theorem \ref{theorem1}. In addition, passivity theory is applied in the quasi-steady state analysis.
%
%At the quasi-steady states, $\eta$ and $w$ are frozen at $\sum_{i=1}^N \varphi_i(y_i)\boldsymbol{1}_N$ and $\bar{w}$ by $\varepsilon=0$, respectively. The reduced system is given by
\begin{equation} \label{rs2}
\begin{aligned}
\dot{\tilde{x}}&=H\tilde{x}+(\bar{B}\otimes I_N)(y-y^*-(F(y)-F(y^*))-\mathbf{A}^T\tilde{\boldsymbol{\mu}}),\\
\dot{\tilde{\boldsymbol{\mu}}}&=-\alpha (L\otimes I_l)\tilde{\boldsymbol{\mu}}-\tilde{z}+\mathbf{A}(y-y^*),\\
\dot{\tilde{z}}&=\alpha (L\otimes I_l)\tilde{\boldsymbol{\mu}}.
\end{aligned}
\end{equation}
For $\tilde{x}$-subsystem, because $H$ is Hurwitz, all eigenvalues of the subsystem have negative real parts. By Lemma 6.1 in \cite{Khalil.2002}, $\tilde{x}$-subsystem is strictly positive real. Furthermore, it follows from Lemma 6.4 in \cite{Khalil.2002} that $\tilde{x}$-subsystem is strictly passive. Thus, there exists a storage function $V_1=\frac{1}{2}\tilde{x}^TP\tilde{x}$ with $P=\blk\{P_1,\ldots,P_N\}$ satisfying that $P_i=P_i^T>0$,  $P_i^TH_i+H_i^TP_i=-Q_i^TQ_i-\epsilon_i P_i$ and $P_i\bar{B}=\bar{C},\forall i\in\mathcal{I}$, where $\epsilon_i$ is a positive constant depending on the largest eigenvalue of $H_i$. Let $\epsilon=\min\{\epsilon_1,\ldots,\epsilon_N\}$. It yields that
\begin{equation}
\begin{split}
\dot{V_1}&\leq (y-y^*-(F(y)-F(y^*))-\mathbf{A}^T\tilde{\boldsymbol{\mu}})^T(\bar{C}\otimes I_N)\tilde{x}\\
&\ \ \  -\frac{1}{2}\epsilon \tilde{x}^TP\tilde{x}.
\end{split}
\end{equation}

For the reduced system \eqref{rs2}, we consider the candidate Lyapunov function $V=V_1+\frac{1}{2}(\|\tilde{\boldsymbol{\mu}}\|^2+\|\tilde{\boldsymbol{\mu}}+\tilde{z}\|^2)$. The derivative of $V$ along the trajectories of the reduced system \eqref{rs2} is
\allowdisplaybreaks[4]
\begin{align*}
\dot{V}&=\dot{V}_1+\tilde{\boldsymbol{\mu}}^T\dot{\tilde{\boldsymbol{\mu}}}+(\tilde{\boldsymbol{\mu}}+\tilde{z})^T(\dot{\tilde{\boldsymbol{\mu}}}+\dot{\tilde{z}})\\
&\leq (y-y^*-(F(y)-F(y^*))-\mathbf{A}^T\tilde{\boldsymbol{\mu}})^T(\bar{C}\otimes I_N)\tilde{x} \\
&\ \ \ -\frac{1}{2}\epsilon \tilde{x}^TP\tilde{x}-\alpha \lambda_2\|\tilde{\boldsymbol{\mu}}\|^2-\|\tilde{z}\|^2\\
&\ \ \ -2\tilde{\boldsymbol{\mu}}^T\tilde{z}+2\tilde{\boldsymbol{\mu}}^T\mathbf{A}(y-y^*)+\tilde{z}^T
\mathbf{A}(y-y^*)\\
&\leq-(\frac{1}{2}\lambda_{min}(P)-\frac{3}{2})\|\tilde{x}\|^2-(\omega-\frac{\|\mathbf{A}\|^2}{2})\|y-y^*\|^2\\
&\ \ \ -(\alpha \lambda_2-2-\frac{\|\mathbf{A}\|^2}{2})\|\tilde{\boldsymbol{\mu}}\|^2-\frac{1}{2}\|\tilde{z}\|^2,
\end{align*}
where $\lambda_{min}(P)$ is the smallest eigenvalue of $P$. Recall the design of matrix $K_i$, it derives that $\epsilon=1$. Next, we discuss the cases of $\omega\geq\frac{\|\mathbf{A}\|^2}{2}$ and $\omega<\frac{\|\mathbf{A}\|^2}{2}$, respectively.

1) In the case of $\omega\geq\frac{\|\mathbf{A}\|^2}{2}$, if $\lambda_{min}(P)>3$ and $\alpha \lambda_2>2+\frac{\|\mathbf{A}\|^2}{2}$, it follows that $\dot{V}<0$, which implies that the states of system \eqref{rs2} can exponentially converge to the origin.

2) In the other case, the sufficient conditions for $\dot{V}<0$ are that $\lambda_{min}(P)>3+2\|\mathbf{A}\|^2$ and $\alpha \lambda_2>2+\frac{\|\mathbf{A}\|^2}{2}$, which ensure that the states of system \eqref{rs2} can exponentially converge to the origin.

The boundary-layer analysis is similar to that in Theorem \ref{theorem1}. Thus, we conclude that there exists a positive constant $\varepsilon^*$ such that for all $\varepsilon \in(0,\varepsilon^*)$, the strategies of all players with strategy-updating rule \eqref{sup2} exponentially converge to the GNE of aggregative game  $G=(\mathcal{I},\Omega,J)$.  $\hfill\blacksquare$

\begin{remark}
In recent work \cite{Bianchi.2019}, the games of multi-integrator agents with coupled constraints were considered but the proposed algorithm was hard to be applied directly here because our setups are different from those in \cite{Bianchi.2019} in three aspects. First, the proposed algorithms 3 and 4 for aggregative games in \cite{Bianchi.2019} are not suitable for agents with inherent complex dynamics or for nonlinear aggregation functions, and these are considered in this note. Second, the algorithm 5 designed in \cite{Bianchi.2019} only considers the coupled constraints, while in this note both local and coupled constraints are considered. In this case, it is much harder to be analyzed. Third, the communication graphs in \cite{Bianchi.2019} are assumed to be undirected and connected, while the graphs in this note are assumed to be directed, strongly connected and weight-balanced. Our assumption is more general.
\end{remark}

\section{Simulations}
In this section, two examples of networked multi-integrator systems are given respectively.

\begin{figure}[!t]\centering
\centering
\includegraphics[width=7cm]{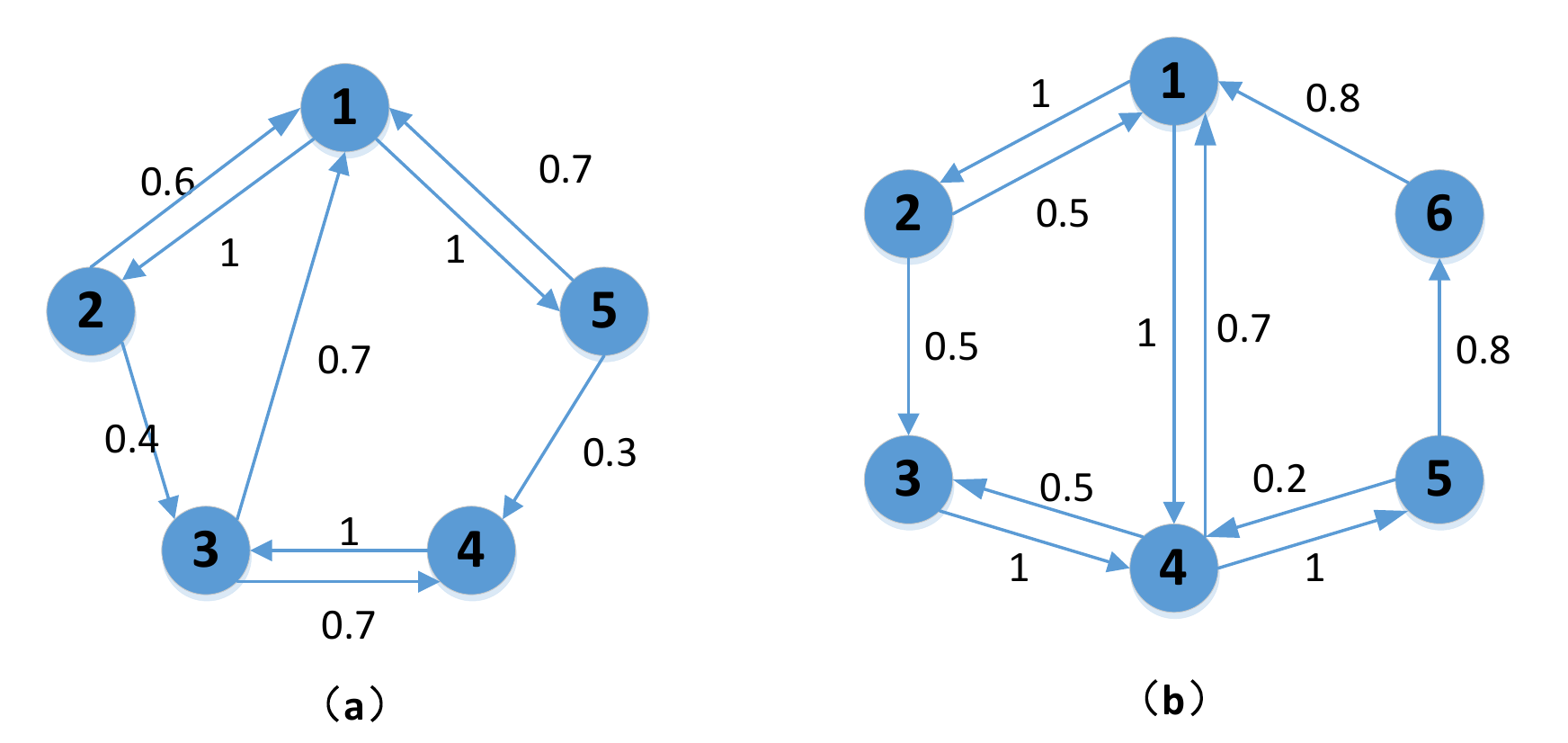}
\caption{Communication graphs}
\label{fig.1}
\end{figure}

\subsection{Formation for Multiple Euler-Lagrange Systems}
Multi-agent formation control can be modeled by a noncooperative game \cite{Stankovic.2012,Lin.2019}. Multiple Euler-Lagrange systems are considered here as an example. The systems with known nonlinearities can be transformed into double-integrator agents. Consider an aggregative game with five Euler-Lagrange systems whose communication graph is depicted in Fig.\ref{fig.1} (a). The cost function of Euler-Lagrange system $i$, $i\in\{1,\ldots,5\}$ is given by
\begin{equation*}
J_i(y_i,y_{-i})=\|y_i-Q_i\|^2+\beta\sum_{i=1}^5\sum_{j=1}^5l_{ij}(y_i-h_i)^T(y_j-h_j),
\end{equation*}
where $y_i\in \mathbb{R}^2$ is the location information of the system $i$'s output, $Q_i$ is a landmark of system $i$, $l_{ij}$ is an element of Laplacian matrix $L$, $h=\col(h_1,\ldots,h_5)$ denotes the desired formation, $\beta$ is a positive constant. The coupled constraint is $ (L\otimes I_2)(y-h)=0$ and the location constraints are $-10\leq y_{i1},y_{i2}\leq 10$.

\begin{figure}[!t]\centering
\centering
\includegraphics[width=6.5cm]{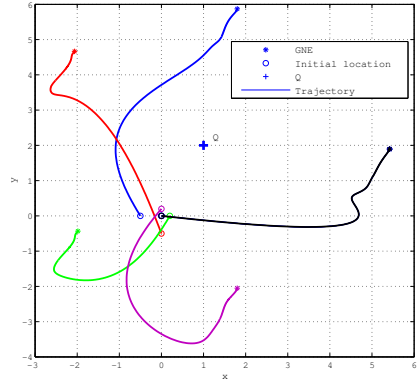}
\caption{Formation of five Euler-Lagrange systems}
\label{fig.2}
\end{figure}

In this case, $Q_1=\cdots=Q_5=[1,2]^T$, the initial location is $y(0)$ $=$ $[0,0; -0.5,0; 0,-0.5; 0.2,0; 0,0.2]^T$, $h_i=[5\cos(2\pi/5(i-1)),5\sin(2\pi/5(i-1))]^T$ and $\beta=5$. The parameters in \eqref{sup1} and \eqref{lme} are selected as $k_i=5$, $\alpha=1$ and $\varepsilon=0.1$ for $i\in \{1,\ldots,5\}$. The trajectories of the five agents are shown in Fig.~\ref{fig.2}. Finally, the five agents form a formation around point $Q$.

\begin{table}
\renewcommand{\arraystretch}{1.3}
\caption{System parameters}
\scalebox{0.7}{
\label{Table 1}
\centering
\begin{tabular}{cccccccccccccc}
\hline\hline \\[-3mm]
~ & $T_{mi}$ & $T_{ei}$ & $K_{mi}(K_{ei})$  & $D_i$ & $H_i$ & $R_i$ & $\alpha_i$ & $\beta_i$ & $\xi_i$ & $P_i(0)$ & $d_i$ & $\Omega_i$\\
\hline
Generator $\sharp$1 & 0.35 & 0.10 & 1.0 & 5.0& 4.0& 0.05& 5 & 12 & 1.0 & 30 & 30 & [20 30]\\
Generator $\sharp$2 & 0.30 & 0.12 & 1.1 & 4.0& 3.5& 0.04& 8 & 10 & 0.5 & 35 & 45 & [45 50]\\
Generator $\sharp$3 & 0.28 & 0.08 & 0.9 & 3.0& 2.8& 0.03& 6 & 11 & 0.8 & 20 & 28 & [25 35]\\
Generator $\sharp$4 & 0.40 & 0.11 & 1.2 & 4.5& 4.2& 0.06& 9 & 11 & 0.7 & 35 & 40 & [30 40]\\
Generator $\sharp$5 & 0.43 & 0.90 & 0.8 & 3.5& 3.0& 0.04& 7 & 13 & 1.1 & 22 & 33 & [20 30]\\
Generator $\sharp$6 & 0.35 & 0.10 & 1.0 & 5.0& 4.0& 0.05& 8 & 14 & 0.6 & 28 & 25 & [20 37]\\
\hline\hline
\end{tabular}}
\end{table}

\subsection{Demand Response of DERs}
In the economic dispatch of power systems, power plants on the supply side compete with each other for the minimum costs, which can be described as an aggregative game. We consider a network of six generation systems communicating with each other over a strongly connected and weighted-balanced digraph depicted in Fig.~\ref{fig.1} (b). The cost function of the generation system $i$ is described by
\begin{equation*}
\begin{split}
J_i(y_i,y_{-i})=c_i(y_i)-p(\sigma)y_i
               =\alpha_i+\beta_iy_i+\xi_iy_i^2-(p_0-a\sigma)y_i,
\end{split}
\end{equation*}
where $y_i\in\mathbb{R}$ is the strategy of the generation system $i$, $c_i(y_i)$ is the generation cost, $p(\sigma)$ is the electricity price, $\alpha_i,\beta_i$, and $\xi_i$ are characteristics of the generation system $i$, $p_0$ and $a$ are constants, and $\sigma=\sum_{i=1}^N y_i$ denotes the linear aggregator.

\begin{figure}[!t]\centering
\centering
\includegraphics[width=6.5cm]{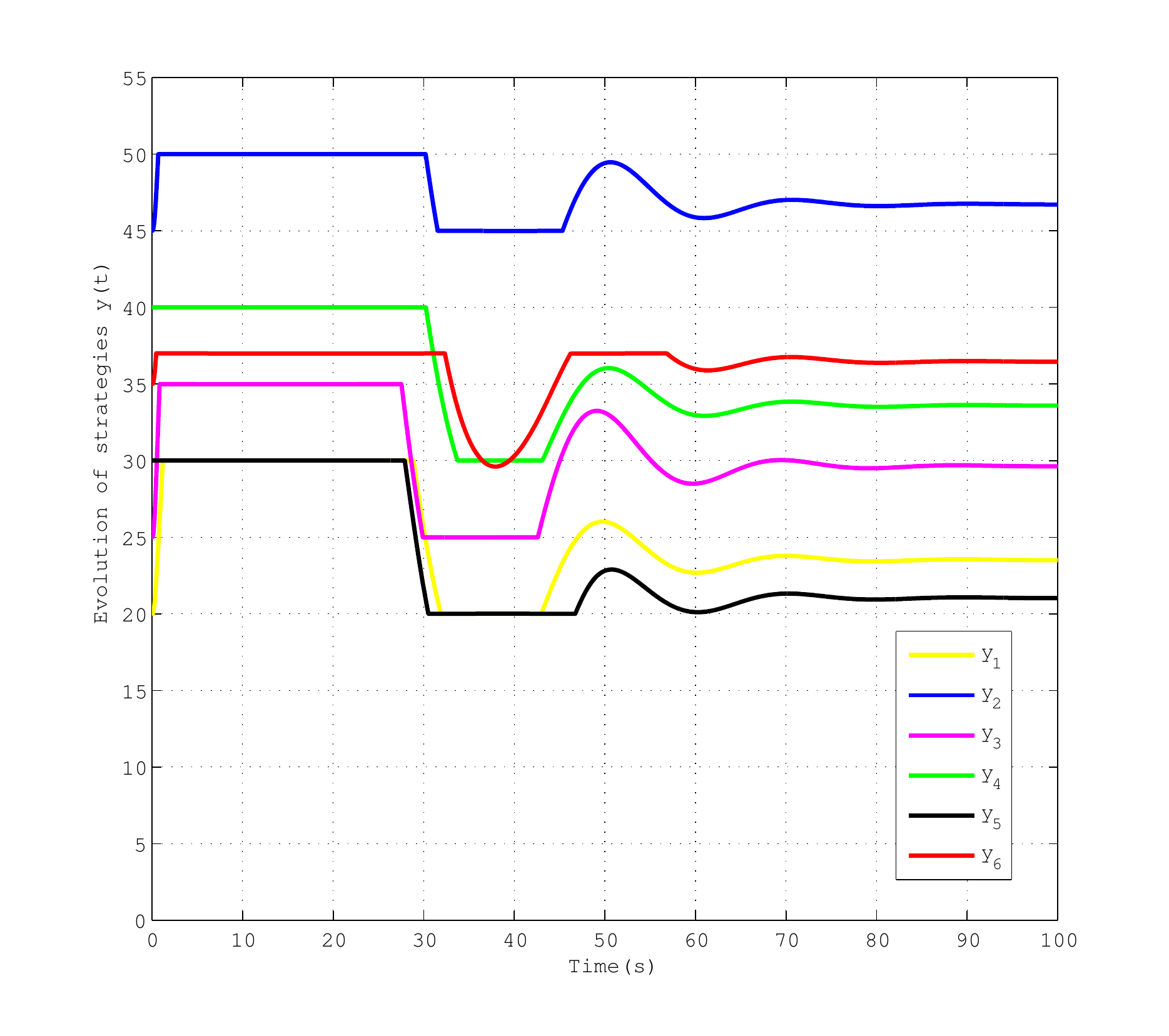}
\caption{Evolution of strategies of multi-integrator agents}
\label{fig.3}
\end{figure}

\begin{figure}[!t]\centering
\centering
\includegraphics[width=6.5cm]{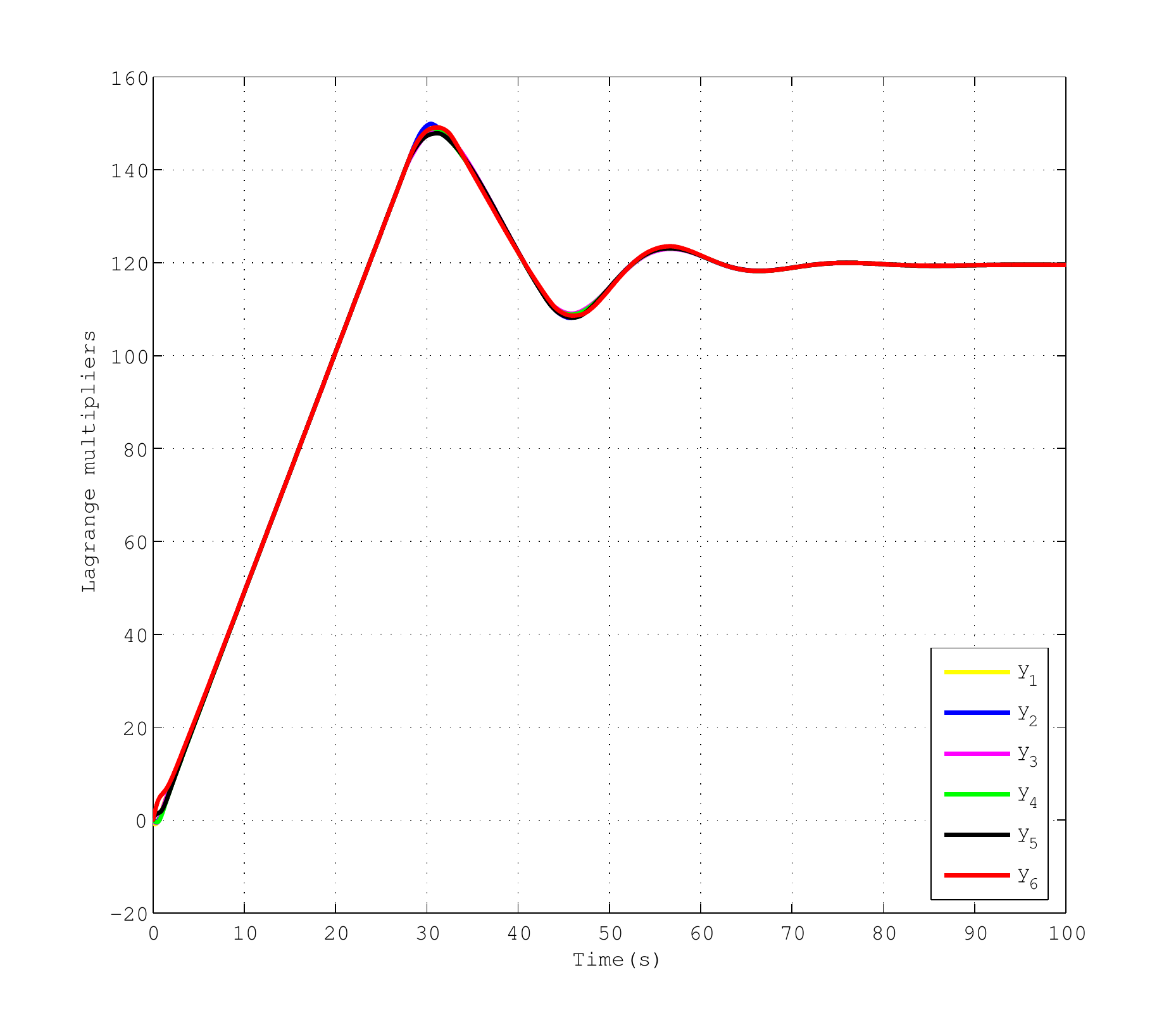}
\caption{Evolution of the Lagrange multipliers of multi-integrator agents}
\label{fig.4}
\end{figure}

The dynamics of the $i$th turbine-generator system are given by (refer to\cite{YiGuo.2000})
\begin{equation*}
\begin{aligned}
\dot{P}_i&=-\frac{1}{T_{mi}}P_i+\frac{K_{mi}}{T_{mi}}X_{ei},\\
\dot{X}_{ei}&=-\frac{K_{ei}}{T_{ei}R_iw_0}w_i-\frac{1}{T_{ei}}X_{ei}+\frac{1}{T_{ei}}u_i,\\
\dot{w}_i&=-\frac{D_i}{2H_i}w_i+\frac{w_0}{2H_i}(P_i-d_i),\\
y_i&=P_{\Omega_i}(P_i),
\end{aligned}
\end{equation*}
where $P_i, X_{ei}$, and $w_i\in\mathbb{R}$ are the power, valve opening and relative speed of the generation system $i$, respectively. $d_i$ is the electricity demand. $T_{mi}$ and $K_{mi}$ are the time constants and the gain of the machine's turbine, respectively. $T_{ei}$ and $K_{ei}$ are the time constants and the gain of speed governor, respectively. $R_i$ is the regulation constant of machine's turbine, $D_i$ is the unit damping constant, $H_i$ is the inertia constant, and $w_0$ is the synchronous machine speed. The parameters of six generation systems are shown in Table \ref{Table 1}. The coupled constraint is $\sum_{i=1}^Ny_i=\sum_{i=1}^Nd_i$, where $[d_1,\ldots,d_6]=[30,45,28,40,23,25]$.

The evolution of strategies and Lagrange multipliers of all players are depicted in Figs.~\ref{fig.3} and \ref{fig.4}, respectively, which illustrate that the strategies of all players evolute in their local constraint sets and converge to the GNE of aggregative game  $G=(\mathcal{I},\Omega,J)$.

\section{Conclusions}
In this paper, we consider aggregative games with coupled constraints in the framework of multi-agent systems, where each agent is described as a multi-integrator system. A distributed strategy-updating rule with two time scales was proposed on the basis of only local information in a weight-balanced digraph. The rule was combined with the coordination of Lagrange multipliers and the estimation of an aggregator, and it deals with coupled constraints and serves the purpose of distributed updating of strategies, respectively.  Via Lyapunov stability theory and singular perturbation theory, the strategies of all players are shown to evolve to the GNE of aggregative game, which are further verified by simulation examples. Future works will focus on  communication costs, social optimal solutions among multiple equilibria, and the influence of stubborn players or cheaters on the games.

\ifCLASSOPTIONcaptionsoff
  \newpage
\fi

\bibliographystyle{IEEEtran}



%\bibliographystyle{IEEEtran}

\end{document}